\documentclass[a4paper, 12pt]{elsart}

\usepackage{amsmath,amssymb}
\usepackage{url}
\usepackage{graphicx}
\usepackage{subfigure}
\usepackage{epsfig}
\graphicspath{{./paper_figures/}}

\def\Cset{\mathbb{C}}

\def\Rset{\mathbb{R}}

\def\cI{{\mathcal I}}

\def\btau{\bar{\tau}}

\def\Tr{{\mathrm {Tr}}}

\newcommand{\eqd}{\mbox{$\; \stackrel{\mbox{{\small \rm def}}}{=} \;$} }

\renewcommand{\Re}{\mathop{\mathrm{Re}}}

\newcommand{\EE}{{\mathcal E}}

\newcommand{\LL}{{\mathcal L}}

\newcommand{\NN}{{\mathcal N}}
\newcommand{\DD}{{\mathcal D}}

\newcommand{\MM}{{\mathcal M}}
\newcommand{\II}{{\mathcal I}}

\newcommand{\Dx}{{\partial _x}}
\newcommand{\Dy}{{\partial _y}}
\newcommand{\DDy}{{\partial^2 _y }}
\newcommand{\lapl}{{\Delta}}
\newcommand{\Ig}{{\int_G}}
\newcommand{\SH}[1]{({#1} )_{H^1}} 
\newcommand{\NL}[1]{\| {#1}\|_{L^2}} 
\newcommand{\SL}[1]{({#1} )_{L^2}} 
\newcommand{\AV}[1]{{\langle {#1} \rangle}} 
\newcommand{\Ho}{{H^1_\perp}} 
\newcommand{\bv}{{\bar v}} 

\newcommand{\ba}[1]{\begin{array}{#1}}
\newcommand{\ea}{\end{array}}
\newcommand{\matrice}[2]{\left[\ba{#1} #2 \ea\right]}

\def\proof{\noindent{\it Proof: }}
\def\endproof{\hspace*{\fill}~$\blacksquare$\par\endtrivlist\unskip}

\newtheorem{theorem}{Theorem}
\newtheorem{lemma}{Lemma}

\newtheorem{definition}{Definition}
\newtheorem{remark}{Remark}

\begin{document}
\begin{frontmatter}
\title{Average consensus problems in networks of agents with delayed communications}
\author[inria]{Pierre-Alexandre Bliman\thanksref{corr}}
\author[inria]{Giancarlo Ferrari-Trecate }
\address[inria]{INRIA, Domaine de Voluceau \\
                Rocquencourt - B.P.105, 78153, Le Chesnay Cedex, France \\
                e-mail: \tt
                \{Pierre-Alexandre.Bliman,Giancarlo.Ferrari-Trecate\}@inria.fr \\
                }
\thanks[corr]{Corresponding author. Tel: +33 1 39 63 55 68, Fax:+33 1 39 63 57 86.}
\begin{abstract}The present paper is devoted to the study of average consensus problems
for undirected networks of dynamic agents having communication delays.
The accent is put here on the study of the time-delays influence: both constant and time-varying delays are considered, as well as uniform and non uniform repartitions of the delays in the network.
The main results provide sufficient conditions (also necessary in most cases) for existence of average consensus under bounded, but otherwise unknown, communication delays.
Simulations are provided that show adequation with these results.
\end{abstract}

\end{frontmatter}
\section{Introduction}
\label{sec:intro}
In the last few years, the study of multi-agent systems has received a
major attention within the control community. Driving applications
include  unmanned aerial vehicles, satellite clusters,
automated highways and mobile robots. In all cases the aim is to
control a group of agents connected through a wireless network. More
precisely, rather than stabilizing the movement of each agent around a
given set point, the goal is to understand how to make the agents
coordinate and self-organize in moving formations. This problem
becomes even more challenging under partial communication protocols,
i.e. when each agent exchanges information only with few others. 

Many works in the literature focused on conditions for guaranteeing
that the agents asymptotically reach a \emph{consensus}, i.e. they agree upon
a common value of a quantity of interest \cite{JLM03}, \cite{TJP03a}, \cite{TJP03b},
\cite{OSM04}, \cite{OSM03b}, \cite{FTBG04}. As an example, in a network of moving vehicles
a form of consensus is represented by alignment, that happens when all vehicles
asymptotically move with the same velocity. In the aforementioned papers, 
consensus problems have been studied under a
variety of assumptions on the network topology (fixed/switching), the
communication protocol (bidirectional or not), additional performance requirements
(e.g. collision avoidance, obstacle avoidance, cohesion), and the control scheme adopted
(also termed \emph{consensus protocol}).
So far, just few works considered consensus problems when communication is affected by
time-delays. Some results for discrete-time agent models are given in \cite{FPP04} and
\cite{AB04}. Two different consensus protocols for continuous-time agent dynamics
have been investigated in  \cite{Mo04} and
\cite{OSM04}. More specifically, assuming that agents behave like integrators
and that communication delays are constant in time and uniform (i.e. they
have the same value in all channels), an analysis of the maximal 
delay that can be tolerated without compromising consensus has been performed in
\cite{Mo04} and \cite{OSM04}. In particular, the protocol adopted
in \cite{OSM04} is capable to guarantee \emph{average consensus} (i.e. the state of each
agent converges, asymptotically, to the average of the initial agent states rather than to an
arbitrary constant) and the authors provide an explicit formula for the largest
transmission delay. 

In the present work we generalize the results of \cite{OSM04} in various
ways. First, we consider uniform and unknown \emph{time-varying} delays and provide upper
bounds to the maximal delay that does not prevent from achieving
average consensus. Second, we derive similar conditions for  networks affected
by \emph{non uniform}, constant or time-varying delays. In the case of non uniform and constant delays,
we also show that if the communication delay between two agents is equal to zero, then
average consensus may achieved irrespectively of the magnitude of all others delays.

The network of agents is modeled in the framework of Partial difference Equations
(PdEs) introduced in \cite{BM04} and used in \cite{FTBG04} analyzing the property of
various linear and nonlinear consensus protocols. PdEs are models that mimic Partial
Differential Equations (PDEs) and provide a mathematical description of the agents
network where ``spatial'' interactions (due to the network structure) and ``temporal''
ones are kept separated and described by operators acting either on space or time. 
Section \ref{sec:tools} provides an introduction to PdEs.
The
main results are presented through Sections \ref{sec:model}-\ref{sec:indep} 
and three simulation experiments are discussed in Section \ref{sec:examples}.
The upper bounds to the maximal tolerable delay depend on some eigenvalues of
suitably defined operators. Although their numerical computation is easy for a given
network, in Section~\ref{se5} we provide their explicit form as a function of the
network size for fully connected and loop-shaped
network.
\section{Tools for functions on graphs}
\label{sec:tools}

The communication network is modeled through an undirected weighted graph $G$ 
defined by a set $\NN=\{1,2,\ldots,N\}$ \emph{nodes} and a
set $\EE\subset \NN\times\NN $ of \emph{edges}. Each node represents
an agent and an edge $(x,y)$ means that the agents $x$ and $y$
share
the information about their states. Agents linked by an arc are
called \emph{neighbors}.
The neighboring relation is
denoted with $x\sim y$ and we assume that $x\sim x$ always
holds. Two nodes $x$ and $y$ are connected by a
\emph{path} if there is a finite sequence $x_0=x, x_1,\ldots, x_n=y$ such that
$x_{i-1}\sim x_{i}$. The graph $G$ is \emph{connected} when each pair of nodes
$(x,y)\in G\times G$ is connected by a path and \emph{complete} if $\EE=G\times G$.

Weights on the communication
links are defined by a function $\omega:\NN\times\NN\mapsto \Rset^+$ with the properties
\begin{subequations}
\begin{align}
  \label{eq:w_sym}
  \omega(x,y)&=\omega(y,x) \\
\omega(x,y)&>0 \Leftrightarrow x\sim y \label{eq:w_pos}
\end{align}
\end{subequations}

Time-varying delays in communications, are elements of the set
\begin{equation}
  \label{eq:dset}
  \DD=\{\tau_i(\cdot):i\in\II\},\quad  \II=\{1,\ldots,r\},\quad  r\leq \frac{N(N-1)}{2}
\end{equation}
where 
$\tau_i:\Rset_+\mapsto \Rset_+$ are piecewise continuous functions. 
A delay is associated to each
edge through the onto function
$T:\EE\mapsto\DD$ verifying $T(x,y)=T(y,x)$. 
The last equality amounts to consider  delays that
are symmetric,  i.e.
 the lags in transmission from $x$ to $y$ and from $y$ to $x$
 do coincide. 
This also motivates the  bound $r\leq \frac{N(N-1)}{2}$ in \eqref{eq:dset}.

\begin{figure}[t]
     \centerline{
     \begin{tabular}[t]{ccc}
       \subfigure[Graph $G$.]{
       \includegraphics[width=4cm]{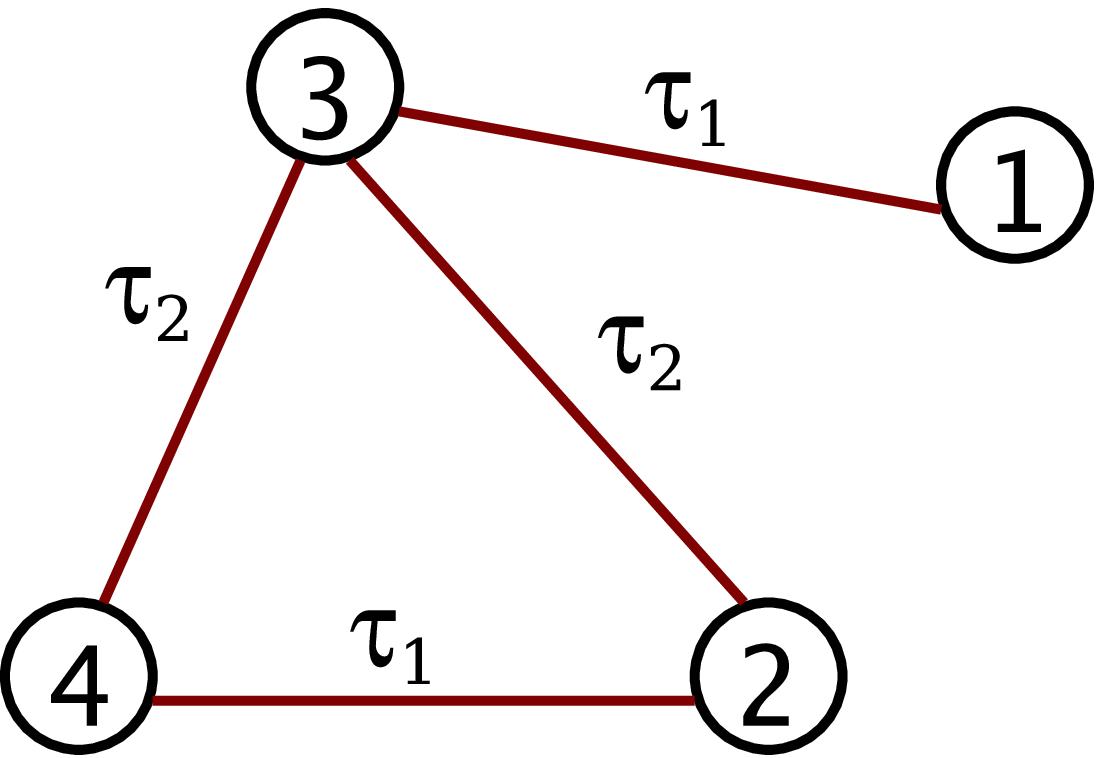}\label{fig:G-a}}&
       \subfigure[Subgraph $G_1$.]{
       \includegraphics[width=4cm]{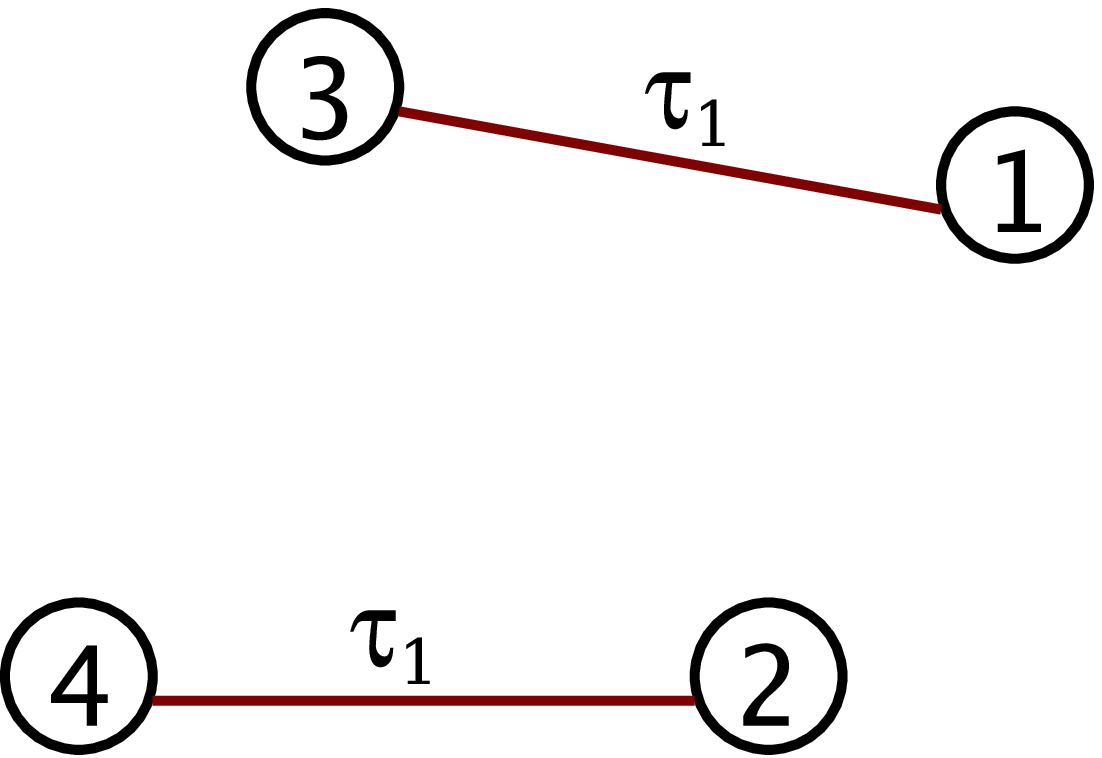}\label{fig:G-b}}&
       \subfigure[Subgraph $G_2$.]{
       \includegraphics[width=4cm]{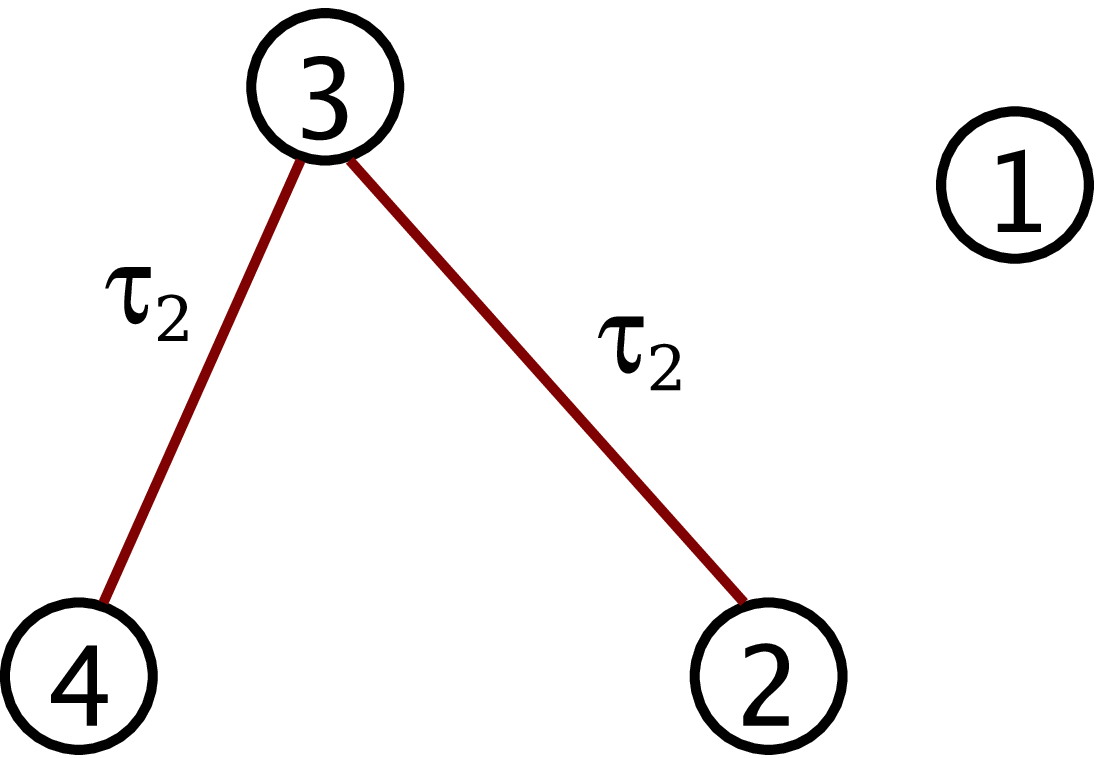}\label{fig:G-c}}
     \end{tabular}
     }
     \caption{\protect{\small A graph and the subgraphs associated to
         delays $\tau_1$ and $\tau_2$}}
     \label{fig:G}
  \end{figure}

Agents linked with the same delay $\tau_i(\cdot)$, define
a subgraph $G_i=(\NN, T^{-1}(\tau_i))$ with associated weights
\begin{equation}
  \label{eq:wi}
  \omega_i(x,y)=
\begin{cases}
\omega(x,y)  & \mbox{if } (x,y)\in T^{-1}(\tau_i) \\
0 & \mbox{otherwise}
\end{cases}
\end{equation}
An example is reported in Figure
\ref{fig:G}. We highlight that the subgraphs $G_i$ may be disconnected
even if $G$ is connected. Moreover, as shown in Figure~\ref{fig:G-c},
some nodes can be isolated. 

We consider now vector
functions $f:\NN\mapsto \Rset^d$ defined over a graph $G$. For instance, $f(x)$
may represent the position or the velocity of the agent $x$ at a fixed
time-instant. Following \cite{BM04}, the \emph{partial derivative} of $f$ is defined as
\begin{equation}
  \label{eq:dxy}
 \Dy f(x)\doteq f(y)-f(x) 
\end{equation}
and enjoys the following  basic  properties:
\begin{subequations}
\begin{align}
\Dy f(x)&=-\Dx f(y) \label{eq:dy-a} \\
\Dx f(x)&=0 \label{eq:dy-b} \\
\DDy f(x)&=\Dy f(y)-\Dy f(x)=-\Dy f(x). \label{eq:dy-c}
\end{align}
\label{eq:dys}
\end{subequations}
The integral and average of $f$ are defined, respectively, as
\begin{equation}
  \label{eq:intG}
  \Ig f ~dx\doteq \sum_{x\in \NN}f(x), \quad  \AV{f}\doteq \frac{1}{N}\Ig f~dx.
\end{equation}
Note that, in \eqref{eq:intG}, ``$dx$'' just indicates the integration variable. The Laplacian of $f$ is given by
\begin{equation}
  \label{eq:Laplacian}
  \lapl f(x)\doteq-\sum_{y\sim x}\omega(x,y)\DDy f(x)=+\sum_{y\sim x}\omega(x,y)\Dy f(x).
\end{equation}
where the last identity follows from \eqref{eq:dy-c}.
In an equivalent way, the Laplacian can be written as
\begin{equation}
  \label{eq:int_Laplacian}
  \lapl f(x)=\int_G\omega(x,y)\Dy f(x)~dy
\end{equation}
The Laplacian operator associated to a  subgraph $G_i$ is
\begin{equation}
  \label{eq:lapl-sub}
 \lapl_i f(x)\doteq\int_{G}\omega_i(x,y)\Dy f(x)~dy
\end{equation}
Since the sets of edges $\{T^{-1}(\tau_i)\}_{i\in\II}$ are a partition of $\EE$, it is
immediate to verify that
\begin{equation}
\label{eq:lapl-split}
\omega(x,y)=\sum_{i\in\II}\omega_i(x,y)\quad\mbox{and}\quad\lapl f= \sum_{i\in\II}\lapl_i f
\end{equation}
In the sequel we summarize the main properties of the Laplacian
operator stated in \cite{BM04}. The driving idea is to mimic functional
analysis tools for studying the classic Laplacian defined on Sobolev
spaces, (see \cite{BM04} and \cite{FTBG04} for further details).
 
We denote with $L^2(G|\Rset^d)$  the Hilbert space composed by all  functions
$f:\NN\mapsto  \Rset^d$ equipped with the scalar product and the norm
\begin{equation}
  \label{eq:L2}
  \SL{f,g}=\Ig f\cdot g,\quad \NL{f}^2=\Ig \|f\|^2
\end{equation}
where $\cdot$ and $\|\cdot\| $ represent the scalar product and the euclidean norm on
$\Rset^d$, respectively. 
Let $H^1(G|\Rset^d)$ be the space collecting all functions in
$L^2(G|\Rset^d)$ with
zero average. We will use the shorthand notation
$L^2$ and $H^1$ when
there is no ambiguity on the underlying domain and range  of the
functions.
If $G$ is connected, $H^1$
is an Hilbert space \cite{BM04} endowed with scalar product
\begin{equation}
  \label{eq:H1-scalar}
  \SH{f,g}=\int_G\int_G\omega(x,y)\Dy f(x)\cdot \Dy g(x)~dxdy. \\
\end{equation}
Apparently, $\Ho$ is the space of \emph{constant}
functions on $G$ and
$\dim(\Ho)=d$. Moreover, the
decomposition $L^2=H^1\oplus \Ho$ is direct. The $L^2$ orthogonal
projection operators on $H^1$ and $\Ho$ will be denoted as $P_{H^1}$
and $P_{\Ho}$, respectively.

The eigenstructure of the Laplacian is completely characterized by
the next Theorem, proved in~\cite{BM04}. 
\begin{theorem}
\label{th:eiglapl}
Let $G$ be a connected graph. Then,
\begin{enumerate}
\item the operator $\lapl:H^1\mapsto H^1$ is symmetric, it has $(N-1)d$ strictly negative
  eigenvalues\footnote{Such eigenvalues will be termed ``the
    eigenvalues of $\lapl$ on $H^1$''.} and the corresponding eigenfunctions form a basis for
  $H^1$;
\item for $f\in L^2$, $\lapl f=0$  if and only if $f\in \Ho$.
\end{enumerate}
\end{theorem}
Theorem~\ref{th:eiglapl} highlights that the Laplacian is invertible
on the subspace $H^1$. Note that when $\lapl$ is defined on $L^2$, it has
$Nd$ eigenvalues. In particular, in view of the decomposition
$L^2=H^1\oplus \Ho$, $(N-1)d$ eigenvalues are those considered
in  point (1) of Theorem~\ref{th:eiglapl} and  the remaining $d$ are
zeros (this property follows directly from  point (2) of
Theorem~\ref{th:eiglapl}).

The next theorem characterizes the eigenvalues of the operators $\lapl_i$.
\begin{theorem}
\label{th:eiglapl-i}
The operators $\lapl_i:H^1\mapsto H^1$, $i\in\II$, are symmetric and negative-semidefinite.
\end{theorem}
\proof
As in the proof of \cite[Lemma 3.1]{BM04}, by direct calculation one
has that $\forall u,v\in L^2$ it holds $\SL{\lapl_i u,v}=-a_i(u,v)$ where $a$ is the
symmetric bilinear form given by
\begin{equation}
  \label{eq:bilin}
  a_i(u,v)\doteq\frac{1}{2}\int_{G_i}\int_{G_i}\omega_i(x,y)\Dy u(x)\cdot \Dy v(x)~dxdy
\end{equation}
Then,  $\SL{\lapl u,v}=- a(u,v)=- a(v,u)=\SL{u,\lapl v}$ that
proves the symmetry of $ \lapl$. For proving that
$\lapl_i$ is negative semidefinite, it is enough to show that $\SL{ - \lapl_i
u,u}\geq 0$, $\forall u\in L^2$. One has
\begin{equation}
  \label{eq:pos}
  \SL{- \lapl_iu,u}=a_i(u,u)=\frac{1}{2}\int_{G_i}\int_{G_i}\omega_i(x,y)\|\Dy u(x)\|^2~dxdy
\end{equation}
where the last term is non negative in view of \eqref{eq:wi} and \eqref{eq:w_pos}.
\endproof

We stress that all the spaces so far introduced are finite
dimensional. This can be seen by noting that
the lifting operator
$\LL:L^2(G|\Rset^d)\mapsto \Rset^{Nd}$ defined as
\begin{equation}
  \label{eq:lift}
  \LL(f)\doteq\matrice{ccc}{f(1)^T &\cdots & f(N)^T}^T
\end{equation}
is an isometry (i.e. bijective and $\NL{f}=\|\LL(f)\|)$ so showing that
$L^2$ is isomorphic to $\Rset^{Nd}$. Roughly speaking, this means that
all concepts introduced in the present section could be
re-written in terms of vector and matrices over $\Rset^{Nd}$.
\begin{definition}
\label{def:mrepr}
  Consider the linear operator $A:L^2(G|\Rset^d)\mapsto L^2(G|\Rset^d)$. Its matrix
  representation is the unique matrix $\MM(A)\in \Rset^{Nd\times Nd}$ that verifies
$\LL(Af)=\MM(A)\LL(f)$, $\forall f\in  L^2(G|\Rset^d)$. 
\end{definition}
The matrix representation of an operator can be used, for instance,
for computing the eigenvalues of $A$ since they coincide with the eigenvalues of
$\MM(A)$, up to their multiplicity.
The operator $\lapl$ is strongly related to the
Laplacian matrix of the graph $G$, defined next (see also \cite{B98}).
In the sequel, the $x$-th row and
the $(x,y)$ element of a matrix $B$ will be denoted with $(B)_x$ and
$(B)_{x,y}$, respectively.
\begin{definition}
\label{def:lapl_matrix}
For a graph $G$, the adjacency matrix $A(G)$ is an $N\times N$ matrix with
entries
\begin{equation}
  \label{eq:adj}
  (A(G))_{x,y}\doteq \begin{cases} \omega(x,y) &\mbox{if } x\sim y
    \mbox{ and } x\neq
  y \\
0 &\mbox{otherwise}\end{cases}
\end{equation}
The valency matrix $V(G)$ is an $N\times N$ diagonal matrix with
entries $(V(G))_{x,x}\doteq\sum_{y\sim x}\omega(x,y)$ and the Laplacian
matrix is $L(G)\doteq A(G)-V(G)$.
\end{definition}
It is easy to verify that 
$\LL(\lapl f(x))=(L(G)\otimes I_d)\LL(f)$, where $\otimes$ is the
Kronecker product and $I_d$ the identity matrix of order $d$ 
Then, $\MM(\lapl)=L(G)\otimes I_d$.

\section{Delayed multi-agent models and PdEs}
\label{sec:model}
Let $v(x,t)\in \Rset^d$ and $u(x,t)\in \Rset^d$, $x\in \NN$, $t\in
\Rset^+$ denote the state and the control input
of agent $x$ at time $t$, respectively. 
When each agent behaves as an integrator,
the collective dynamics is described by the equation $\dot
v(x,t)=u(x,t)$, where the dot operator indicates the time-derivative. In this
paper we consider \emph{delayed Laplacian protocols} of the type
$u=\sum_{i\in\II}\lapl_iv(x,t-\tau_i(t))$ yielding  the collective dynamics
\begin{equation}
  \label{eq:f_dynamics}
  \dot v=\sum_{i\in\II}\lapl_iv(x,t-\tau_i(t))
\end{equation}
Formula (\ref{eq:f_dynamics}) defines a \emph{time-delay Partial
  difference Equation (PdE)} (see \cite{FTBG04} for a general definition of
PdEs) whose solution depends on the initial conditions.
As for linear time-delay systems, if all delays are bounded by a constant
$\bar \tau$, the initial condition may be given in form of a function
$\tilde v(x,t)\in \Rset^d$, $x\in \NN$, $t\in [-\bar \tau,0 ]$, continuous in
$t$.

As shown in  \cite{FTBG04}, PdEs can be always recast into Ordinary
Differential Equations  by using the lifting operator \eqref{eq:lift}.
Then, it is not surprising that linear time-delay PdEs 
inherit all the properties of linear time-delay
systems. As an example, if all delays are constant in time,
the characteristic equation associated to
\eqref{eq:f_dynamics}, is
\begin{equation}
  \label{eq:char}
  E(s)\doteq sI-\sum_{i\in\II}e^{-s\tau_i}\lapl_i=0,\quad s\in \Cset
\end{equation}
where $I$ is the identity operator on $L^2$. Then, many 
properties of the network of agents can be characterized in terms of the \emph{poles} of
\eqref{eq:f_dynamics}, i.e. the roots of \eqref{eq:char}. 
We outline that if the delays are constant, model \eqref{eq:f_dynamics} coincides with
the network dynamics considered in Section 10 of \cite{OSM04}.

The main goal of the present work is to investigate when
\eqref{eq:f_dynamics} guarantees average consensus.
\begin{definition}
  \label{def:align}
The network dynamics achieves average consensus if $v\to\AV{v(\cdot,0)}$ as $t\to+\infty$.
\end{definition}

In absence of delays,  $u$ results in the \emph{Laplacian protocol},   and the
PdE \eqref{eq:f_dynamics} reduces to the \emph{heat equation}
\begin{equation}
  \label{eq:f_heat}
  \dot v=\lapl v\quad v(\cdot,0)=\tilde v\in L^2
\end{equation}
The consensus properties of Laplacian protocols have been analyzed in various works.
In particular,
A. Jadbabaie \emph{et al.}  \cite{JLM03} proved that the Laplacian protocol
is able to guarantee average consensus under various assumption on the network topology. 
A formal
analysis of the PdE \eqref{eq:f_heat} has been carried out in
\cite{FTBG04}, where it has been also shown that the Laplacian
protocol can guarantee consensus even when the agent dynamics are
perturbed by exponentially decreasing errors and/or an agent acts as
the leader of the group.

In order to highlight the rationale we will use for analyzing the PdE \eqref{eq:f_dynamics},
let us summarize
the main results of \cite{FTBG04} for the collective dynamics
\eqref{eq:f_heat}. Decomposing the velocities as $v=v_1+\bv$, $v_1\in
H_1$, $\bv=\AV{v(\cdot,t)}\in\Ho$, one can show, through a simple variational
technique, that the velocity components fulfill the dynamics
\begin{subequations}
\begin{align}
\dot \bv &= 0 \label{eq:bv} \\
\dot v_1&=\lapl v_1 \label{eq:v1}
\end{align}
\label{eq:vel}
\end{subequations}
thus proving that the spaces $H_1$ and $\Ho$ are positively invariant
for \eqref{eq:f_heat}. In particular, equation \eqref{eq:bv}
highlights that the average velocity of the agents is constant in
time. Then an \emph{exponentially stable} average consensus is achieved if the
origin of \eqref{eq:v1} is exponentially stable, a fact that
can be easily shown by  exploiting the characterization of the eigenvalues of $\lapl$ on
$H^1$ given in 
Theorem~\ref{th:eiglapl}.
In \cite{FTBG04} it is also shown that average consensus
can be intuitively expected on the basis of
the physical analogy between
\eqref{eq:f_heat} and the classic heat equation.

For the delayed model \eqref{eq:f_dynamics}, we will adopt a similar
argument. The next Lemma provides the dynamics of the $v_1$ and $\bv$
components.
\begin{lemma}
\label{lem:dec}
The function $v$ is
solutions to the PdEs \eqref{eq:f_dynamics}  if and only if $v_1$ and $\bv$, are solutions to the PdEs
\begin{equation}
\label{eq:PdEs_dec}
 \Sigma_{1}:~ \dot v_1  =\sum_{i\in\II}\lapl_i v_1(x,t-\tau_i(t)), \quad \quad  \bar \Sigma:~ \dot \bv  =0 
\end{equation}
equipped with the initial conditions $v_1(\cdot,t)=P_{H^1}\tilde v(\cdot,t) $,
$\bar v(t)=P_{\Ho}\tilde v(\cdot,t)$ for $t\in[-\bar \tau,0]$.
\end{lemma}
\proof
To prove the result, we use a variational argument by testing each side
of \eqref{eq:f_dynamics} against all $c\in\Ho$. This means that we take
the integrals
\begin{equation}
  \label{eq:test}
  \int_Gc \cdot \dot v  ~dx = \int_G c \cdot \sum_{i\in\II}\lapl_iv(x,t-\tau_i(t)) ~dx
\end{equation}
By using \eqref{eq:lapl-sub}, the right side of \eqref{eq:test} can be written as
$\sum_{i\in\II}S_i$, where
\begin{equation}
  \label{eq:test1}
S_i=\int_G c \cdot\int_G \omega_i(x,y)\Dy v(x,t-\tau_i(t)) ~dydx
\end{equation}
From \eqref{eq:w_sym} and \eqref{eq:dy-a}, the functions
$g_i(x,y)=\omega_i(x,y)\Dy v(x,t-\tau_i(t))$ are antisymmetric,
i.e. $g_i(x,y)=-g_i(y,x)$. Then, each integral $S_i$ can be expanded into
sums containing only terms of the type $c \cdot (g_i(x,y)+g_i(y,x))$ that
are all identically equal to zero. The fact that $ \int_Gc \cdot \dot v  ~dx =0$,
 $\forall c\in \Ho$  corresponds to the condition $P_{\Ho}\dot v=0$,
 or, equivalently, to $\dot \bv=0$, thus obtaining the dynamics $\bar \Sigma$.
From \eqref{eq:f_dynamics} we have
\begin{equation}
  \label{eq:split-v}
  \dot v_1+\dot \bv= \sum_{i\in\II}\lapl_iv_1(x,t-\tau_i(t))+\sum_{i\in\II}\lapl_i\bv(x,t-\tau_i(t))
\end{equation}
and the dynamics $\Sigma_1$ follows from $\dot \bv=0$ and $\lapl_i\bv=0$.
\endproof
Lemma~\ref{lem:dec} shows that the spaces $H^1$ and $\Ho$ are
positively invariant for the PdE \eqref{eq:f_dynamics}. Moreover, as for
\eqref{eq:f_heat}, the average state $\bar v$ is constant in time and equal to
$\AV{\tilde v(\cdot,0)}$.
Then, the problem of checking average consensus is reduced to 
the problem of proving that $v_1\to 0$ as $t\to \infty$. 
We also say that average consensus is 
\emph{globally exponentially or asymptotically stable} if the zero solution
to $\Sigma_1$ enjoys the same property, i.e. it is 
exponentially or asymptotically stable 
for all initial conditions $\tilde v(\cdot,t)\in H^1$, $t\in [-\bar \tau,0 ]$, continuous in
$t$.


For subsequent use, we introduce the operator norm  $\|\Delta\|\doteq \max_{u\in H^1}\frac{(u,\lapl u)}{(u,u)}=|\lambda_{min}|$
where $\lambda_{min}$ is the minimal eigenvalue of the Laplacian on $H^1$. 
Similarly, by recalling that $\lapl$ is
invertible on $H^1$, one has $\|\Delta^{-1}\|^{-1}=|\lambda_{max}|$.

\section{The case of uniform delays}
\label{se2}

In this section, we analyze the stability properties of the dynamics $\Sigma_1$
when the delay is \emph{uniform} in the network, i.e. when  $\cI$ is
a singleton.
We start with the simpler case of time-invariant delays, considered also in
\cite{OSM04}.
The results of the next Theorem coincide with those of Theorem 10 in \cite{OSM04}, but
are proved through a different argument, i.e. the diagonalization of the Laplacian
operator on $H^1$. 
\begin{theorem}[Constant delay]
\label{th2}
The zero function is a globally exponentially  stable solution to the PdE
\begin{equation}
\label{syst}
\dot v_1(x,t) = \Delta v_1(x,t-\tau),\quad v_1\in H^1
\end{equation}
for all possible $\tau\leq\btau$,
{\em if and only if}
\begin{equation}
\label{cons}
\btau < \frac{\pi}{2\|\Delta\|}\ .
\end{equation}
\end{theorem}

\begin{proof}
In view of Theorem~\ref{th:eiglapl}, the Laplacian can be diagonalized
on $H^1$. Let $\{\psi_i\}_{i=1}^{(N-1)d}$ be an orthonormal set of eigenfunctions of
$\lapl$ forming a basis for $H^1$ and associated to the eigenvalues
$\{\lambda_i\}_{i=1}^{(N-1)d}$. Then $v_1(x,t)=\sum_{i=1}^{(N-1)d}\alpha_i(t)\psi_i(x)$ for
suitable functions $\alpha_i:\Rset^+\mapsto\Rset$. By testing each side of
\eqref{syst} against $\psi_j$ we form the integrals
\begin{equation}
  \label{eq:test2}
  \int_G\left(\sum_{i=1}^{(N-1)d}\dot \alpha_i(t)\psi_i(x)\right)\cdot \psi_j(x)~dx=
  \int_G\left(\sum_{i=1}^{(N-1)d}\alpha_i(t-\tau)\lapl \psi_i(x)\right)\cdot \psi_j(x)~dx
\end{equation}
By Theorem \ref{th:eiglapl}, formula \eqref{eq:test2} reduces to 
\begin{equation}
  \label{eq:diag}
  \dot \alpha_j(t)=\lambda_j\alpha_j(t-\tau)
\end{equation}
System \eqref{eq:diag} is a first-order linear time-delay
system. Since $\lambda_j<0$,  according to \cite[Theorem A.5]{HAL93}, system
\eqref{eq:diag} is exponentially
stable if
and only if $\tau<\frac{\pi}{2|\lambda_j|}$. Then, the PdE \eqref{syst} is
exponentially stable if and only if all systems \eqref{eq:diag}, for $j=1,\ldots,(N-1)d$ are
exponentially stable, i.e. if \eqref{cons} holds.
\end{proof}

\begin{remark}
\label{re1}
For $h\geq0$, it may be of interest to quantify the largest delay $\btau_h$ 
for which an exponential decay rate $h$ is guaranteed for the solutions to \eqref{syst}.
By using the diagonalization procedure in the proof of Theorem~\ref{th2},
and noting that systems \eqref{eq:diag} are asymptotically stable for $\tau=0$, the quantity $\btau_h$ is the smallest $\tau>0$ for which
there exists $\omega\in\Rset$ such that
\begin{equation}
  \label{eq:h}
  -h+j\omega -e^{\tau(h-j\omega)}\| \lapl \|=0.
\end{equation}
Using the fact that if $\omega$ verifies \eqref{eq:h} also $-\omega$ verifies \eqref{eq:h}, one obtains
$\omega=\sqrt{\| \lapl \|^2e^{2h\tau}-h^2}$. Furthermore, summing up the real and complex
parts of equation \eqref{eq:h} yields
\begin{equation*}
\btau_h
\doteq \min\left\{
\tau\geq 0\ :\  \|\Delta\|e^{h\tau}\cos\left(\tau\sqrt{\|\lapl\|^2
 e^{2h\tau}-h^2}\right)=h \right\} \ .
\end{equation*}
Note that the map $h\mapsto\btau_h$ is decreasing, with $\btau_0=\pi/2\|\Delta\|$, and $\btau_{\|\Delta\|}=0$.
\end{remark}

We consider now  the case of a single {\em time-varying\/} delay.
\begin{theorem}[Time-varying delay]
\label{th3}
The zero solution is a globally  exponentially  stable solution to the PdE
\begin{equation}
\label{syst2}
\dot{v_1}(x,t) = \Delta v_1(x,t-\tau (t)),\quad v_1(\cdot,t)\in H^1
\end{equation}
for all piecewise continuous delays $\tau(t)$ verifying $0\leq\tau(t)\leq \btau$,
{\em if and only if}
\begin{equation}
\label{vary}
\btau < \frac{3}{2\|\Delta\|}\ .
\end{equation}
\end{theorem}

\begin{proof}
As in the proof of Theorem \ref{th2}, diagonalization of the Laplacian
on $H^1$
leads to the study of the first-order systems
$\dot \alpha_i=\lambda_i\alpha_i(t-\tau(t))$, for any eigenvalue $\lambda_i$ of $\Delta$
on $H^1$.

The conclusion is then deduced from a classical result initially published
in \cite{MYS51} and \cite{YOR70}, (see also \cite[p.\ 164]{HAL93} and the references therein).
\end{proof}

If the nominal collective model is the PdE \eqref{eq:f_heat},
Theorems \ref{th2} and \ref{th3} characterize the robustness of
average consensus with respect to different delay models. In particular, 
the bounds given in Theorems \eqref{cons} and \eqref{vary} do not depend
upon the precise structure of the communication network but only upon
the magnitude of $\|\Delta\|$. In other words, by interpreting $G$ as 
the ``spatial'' domain of the PdEs \eqref{syst}
and \eqref{syst2}, bounds \eqref{cons} and \eqref{vary} relate the
maximal tolerated delays to a spatial feature. 
Explicit formulas for $\|\lapl\|$ in
the case of complete and loop-shaped networks are given in
Sections~\ref{subsec:complete} and \ref{subsec:loop},
respectively.
Other results linking the graph structure with
the eigenvalues of the Laplacian operator 
can be found in \cite{GMS90}, \cite{GM94} and \cite{M94}.

We also outline that the constant in \eqref{vary} is smaller than the corresponding one in
\eqref{cons}, the greater conservativity arising from the time-varying
nature of the delay. However, the bound \eqref{vary} is the best
possible one since the corresponding stability
condition is \emph{necessary and sufficient}.

\section{The case of non-uniform delays}
\label{se3}

In this Section, we generalize the results of Section~\ref{se2}
to the case where the delays do not take a common value in the whole
network.
Let us consider first the case of constant delays. The next Theorem
provides a robust stability result for \emph{all} possible delays
$\tau_i$ within the interval $[0,\bar \tau]$.
Quite remarkably, the bound \eqref{cons} still gives a necessary and
sufficient condition for stability.

\begin{theorem}[Constant delays]
\label{th4}
The zero solution is a globally exponentially  stable solution to the equation
\begin{equation}
\label{equa}
\dot{v_1}(x,t) = \sum_{i\in\cI}\Delta_i v_1(x,t-\tau_i),\quad v_1(\cdot,t)\in H^1
\end{equation}
for all possible $\tau_i\leq\btau$, $i\in\II$, 
{\em if and only if\/} \eqref{cons} holds.
\end{theorem}

\proof
By considering the case $\tau_i=\btau$, $i\in\cI$, Theorem \ref{th2}
shows that the upper bound to the tolerated delay 
cannot be larger than $\pi/2\|\Delta\|$.

We prove by contradiction
that \eqref{cons} implies asymptotic stability.
Assume that \eqref{equa} is not asymptotically stable.
For zero delays, the PdE \eqref{equa} reduces to PdE
\eqref{eq:v1} whose global exponential stability has been proved in
\cite{FTBG04}.
By continuity of the  poles of \eqref{equa} with respect to the
delays, there exists a choice of the $\tau_i\in [0,\btau]$, $i\in\cI$, for which
the PdE \eqref{equa} has a purely imaginary pole $j\omega$, $\omega\in\Rset\setminus\{ 0\}$.
In other words, there exists a nonzero  eigenfunction $u\in H^1(G|\Cset^d)\backslash\{0\} $
 such that
\[
\left(
j\omega I -\sum_{i\in\cI}\Delta_i e^{-j\tau_i\omega}
\right) u=0\ ,
\]
where $I$ is the identity on $H^1$. This implies that
\begin{equation}
\label{part}
j\omega+ \sum_{i\in\cI}\alpha_i e^{-j\tau_i\omega} = 0\ ,
\end{equation}
where, by denoting with $u^*$ the complex conjugate of $u$, one has
\[
\alpha_i \doteq -\frac{\int_Gu^*\cdot \Delta_iu}{\int_Gu^*\cdot u},\ i\in\cI\ .
\]
Notice that, by Theorem~\ref{th:eiglapl-i}, the operators $\Delta_i$ are
symmetric, negative semidefinite on $H^1$. Thus, one has $\alpha_i\geq 0$ for all $i\in\cI$.
On the other hand,
\begin{equation}
\label{ineq}
\sum_{i\in\cI}\alpha_i \leq \|\Delta\|\ .
\end{equation}

Considering the real and imaginary parts of \eqref{part}, we deduce that
\begin{subequations}
\begin{gather}
\label{inea}
\sum_{i\in\cI}\alpha_i\cos \omega\tau_i = 0,\\
\label{ineb}
\omega - \sum_{i\in\cI}\alpha_i\sin \omega\tau_i = 0\ .
\end{gather}
\end{subequations}
From \eqref{ineq} and \eqref{ineb}, one gets that $|\omega|\leq\|\Delta\|$,
whence $|\omega| \tau_i\leq\|\Delta\| \btau<\pi/2$.
In these conditions, the terms $\cos \omega\tau_i$ appearing in \eqref{inea} are
{\em all\/} positive.
Note that not all the coefficients $\alpha_i$ can be zero,
otherwise $\int_Gu^*\cdot \Delta u=0$, which contradicts the fact that $u\in H^1\setminus\{ 0\}$.
Therefore \eqref{inea} is impossible,
and we are thus led to a contradiction.
This proves that if $\tau_i\leq\btau$ and \eqref{cons} holds, then, 
the PdE \eqref{equa} is globally asymptotically stable.
\endproof

We stress once more the robustness flavor of Theorem \ref{th4}, that
requires just the knowledge of a common upper
bound $\bar \tau$  on the (unknown) delays $\tau_i$. On the other hand,
there may exist  combinations of delays $\tau_i$ such that
$\tau_i\geq\bar \tau$, for some $i\in \II$, 
but the PdE \eqref{equa} remains asymptotically stable. An example is
provided in Section~\ref{sec:examples}.

The argument used in the proof of Theorem \ref{th3} 
does not seem to extend to the case of non-stationary delays.
In this case, the next Theorem provides a sufficient stability condition.

\begin{theorem}[Time-varying delays]
\label{th5}
The zero solution is a globally  stable solution to the PdE \eqref{eq:f_dynamics}
for all nonnegative, piecewise continuous delay $\tau_i(t)$ verifying $0\leq\tau_i(t)\leq \btau$,
{\em if}
\begin{equation}
\label{eq:bound-vary-nonuniform}
\btau < \frac{1}{\sum_{i,i'\in\cI}\|\Delta_i\Delta_{i'}\| \|\Delta^{-1}\|}\ .
\end{equation}
\end{theorem}

Theorem \ref{th5} is a direct consequence of the following stronger result.
\begin{lemma}
\label{le1}
The zero solution is a globally stable solution to the equation
\eqref{eq:f_dynamics} for all piecewise continuous nonnegative delay
$\tau_i(t)$ verifying
\begin{equation}
\label{trac}
\sup_{t\in [0,+\infty)}
\left(
\sum_{i,i'\in\cI}\tau_i(t)\|\Delta_i\Delta_{i'}\|
\right)
< 1/\|\Delta^{-1}\|\ .
\end{equation}
\end{lemma}
\proof
One may write
\begin{align}
\dot{v_1} & = \sum_{i\in\cI} \Delta_i v_1(x,t-\tau_i(t))
= \Delta v_1(x,t) +\sum_{i\in\cI} \Delta_i \left(
v_1(x,t-\tau_i(t))-v_1(x,t)
\right)\nonumber \\
& = \Delta v_1(x,t) -\sum_{i\in\cI} \Delta_i
\int_{t-\tau_i(t)}^t \dot{v_1}(x,s)\ ds \nonumber \\
& = \Delta v_1(x,t) -\sum_{i,i'\in\cI} \Delta_i
\int_{t-\tau_i(t)}^t \Delta_{i'}v_1(x,s-\tau_{i'}(s))\ ds\ . \label{eq:dot-v1}
\end{align}

Let $\rho$ be the constant
\[
\rho\eqd 
\sup_{t\in [0,+\infty)}\left(
\sum_{i,i'\in\cI}\tau_i(t)\|\Delta_i\Delta_{i'}\|
\right)
\|\Delta^{-1}\| < 1 \ .
\]
By assumption, there exists a real number $q$ in the non empty interval $(1,1/\rho)$.

We apply now Razumikhin theorem, see e.g.\
\cite[Theorem 4.2]{HAL93} or \cite[Chapter 4]{KOL99}. 
For $V(v_1)\doteq \frac{1}{2}\NL{v_1}$,
assume that $V(v_1(\cdot,s))<q V(v_1(\cdot,t))$ for all $s\in [t-2\btau,t]$.
From the expression of $\dot{v_1}$ in \eqref{eq:dot-v1}, one deduces
that, 
along the trajectories of \eqref{eq:f_dynamics}, it holds
\begin{align*}
\frac{d[V(v_1)]}{dt}
& \leq \left(
-\|\Delta^{-1}\|^{-1} + q \left(
\sum_{i,i'\in\cI}\tau_i(t)\|\Delta_i\Delta_{i'}\|
\right)
\right) \NL{v_1} \\
&\leq
\|\Delta^{-1}\|^{-1}
(-1+q\rho)
\NL{v_1}\ .
\end{align*}
Then, $\frac{d[V(v_1)]}{dt}<0$ and the thesis follows.
\endproof

The results in Theorems \ref{th2}, \ref{th3}, \ref{th4} and \ref{th5} are summarized
in Table \ref{ta1}.

\begin{table}[h]
\begin{center}
\begin{tabular}{c||c|c|}
$\btau$ & Uniform delays & Non-uniform delays \\
\hline
& & \\
Time-invariant delays & $\frac{\pi}{2\|\Delta\|}$\quad {Th.\ \ref{th2}, \sc (e)}
& $\frac{\pi}{2\|\Delta\|}$\quad {Th.\ \ref{th4}, \sc (e)}\\
& & \\
\hline
& & \\
Time-varying delays & $\frac{3}{2\|\Delta\|}$\quad {Th.\ \ref{th3}, \sc (e)}
& $\frac{1}{\sum_{i,i'\in\cI}\|\Delta_i\Delta_{i'}\|\ \|\Delta^{-1}\|}$
\quad {Th.\ \ref{th5}, \sc (s)}\\
& &\\
\hline
\end{tabular}
\caption{Bounds on the worst-case stabilizing delay.
{\sc e}: exact, {\sc s}: sufficient.}
\label{ta1}
\end{center}
\end{table}
\begin{remark}
By comparison with \eqref{vary}, the bound \eqref{eq:bound-vary-nonuniform} depends in a more involved manner upon the
structure of the communication network. 
We also highlight that the bound \eqref{eq:bound-vary-nonuniform} may be bounded from above by the simpler quantity
$1/(\ \Tr \Delta)^2\|\Delta^{-1}\|$, where $\Tr \Delta$ is the trace of the Laplacian
on $H^1$.\\
For checking that the results of Theorems \ref{th3} and \ref{th5} are
coherent, one can use the following inequalities
\begin{align*}
&\sum_{i,i'\in\cI}\|\Delta_i\Delta_{i'}\| \|\Delta^{-1}\|
\geq
\left\|
\sum_{i,i'\in\cI} \Delta_i\Delta_{i'}
\right\| \|\Delta^{-1}\|= \\
&= \left\|
\left(
\sum_{i\in\cI}\Delta_i
\right)^2
\right\| \|\Delta^{-1}\|
= \|\Delta^2\| \|\Delta^{-1}\|
\geq \|\Delta\|\ ,
\end{align*}
that imply \eqref{eq:bound-vary-nonuniform}$\leq
\frac{1}{\|\lapl\|}\leq\frac{3}{2\|\lapl\| }$.
Also, we highlight the trade-off between stability with large delays
on the one hand, and large decay-rate of the solutions on the other hand:
the first one requires a small $\|\Delta\|$, whereas the second one requires a large
$\|\Delta^{-1}\|^{-1}\leq\|\Delta\|$.
\end{remark}

\section{A delay-independent condition for average consensus}
\label{sec:indep}
According to the standard terminology in time-delay systems, all the
results presented in Sections~\ref{se2} and \ref{se3} are
``delay-dependent'' in the sense that they guarantee average consensus
when all the communication delays are upper-bounded by
a suitable value $\bar \tau$. 
Next, we show that if a single delay is zero, average consensus may
be achieved irrespectively of the magnitude of all other delays.
In this sense, we provide a ``delay-independent'' condition for
average consensus.
For two operators $A$ and $B$ from  $L^2$ to  $L^2$ the inequality $A>B$
on $H^1$ means that 
\begin{equation}
\label{eq:ineq}
\forall u\in H^1\setminus\{ 0\},\
\SL{u,(A-B) u}> 0\ .
\end{equation}
\begin{theorem}
Consider the PdE \eqref{equa} and assume that $\tau_{i'}=0$ for an index
$i'\in\cI$.
If 
\begin{equation}
\label{delt1}
\Delta_{i'} < \sum_{i\in\cI\setminus\{ i'\}} \Delta_i \quad \text{ on $H^1$},
\end{equation}
then, the zero solution is a globally exponentially stable solution to \eqref{equa}
for any $\tau_i\geq 0$, $i\in\cI\setminus\{ i'\}$.
Conversely, if the zero solution to system \eqref{equa} with
$\tau_{i'}=0$ is globally asymptotically stable
for any $\tau_i\geq 0$, $i\in\cI\setminus\{ i'\}$, then
\begin{equation}
\label{delt2}
\Delta_{i'} \leq \sum_{i\in\cI\setminus\{ i'\}} \Delta_i \quad \text{ on $H^1$}\ .
\end{equation}
\end{theorem}

\begin{proof}
Assume first that \eqref{delt1} holds.
Then, for any $z_i\in\Cset$, $i\in\cI\setminus\{ i'\}$,
such that $|z_i|\leq 1$, one has
\[
\left(
\Delta_{i'} + \sum_{i\in\cI\setminus\{ i'\}}z_i \Delta_i
\right)
+\left(
\Delta_{i'} + \sum_{i\in\cI\setminus\{ i'\}}z_i \Delta_i
\right)^*
\leq 2 \left(
\Delta_{i'} - \sum_{i\in\cI\setminus\{ i'\}} \Delta_i
\right)
< 0 \ \mbox{ on } H^1\ .
\]
Consequently, $\Delta_{i'} + \sum_{i\in\cI\setminus\{ i'\}}z_i \Delta_i$ has
only eigenvalues with strictly negative real part on $H^1$ for any $|z_i|\leq 1$.
This implies that, for any $s\in\Cset$ with $\Re s\geq 0$,
the following inequality holds on $H^1$:
\begin{multline}
\label{eq:ineq2}
(s+s^*)I \geq 0 >
2\left(\Delta_{i'} + \sum_{i\in\cI\setminus\{ i'\}} \Delta_i\right)\\
\geq \left(\Delta_{i'} + \sum_{i\in\cI\setminus\{ i'\}} z_i\Delta_i\right)+\left(\Delta_{i'} + \sum_{i\in\cI\setminus\{ i'\}} z_i\Delta_i\right)^*\ .
\end{multline}
In particular, one can choose $z_i=e^{-s\tau_i}$ in \eqref{eq:ineq2}, because $|e^{-s\tau_i}|\leq 1$ when $\Re s\geq 0$ and $\tau_i\geq 0$.
Thus,
\begin{equation}
  \label{eq:ineq3}
  \left(sI-\Delta_{i'} - \sum_{i\in\cI\setminus\{ i'\}} e^{-s\tau_i}\Delta_i\right)+\left(sI-\Delta_{i'} - \sum_{i\in\cI\setminus\{ i'\}} e^{-s\tau_i}\Delta_i\right)^*>0
\end{equation}
for any $\tau_i\geq 0$, $i\in\cI\setminus\{ i'\}$.
As a consequence, all the roots of the characteristic equation $\det (sI-\Delta_{i^*} - \sum_{i\in\cI\setminus\{ i^*\}} e^{-s\tau_i} \Delta_i)= 0$ of \eqref{equa} have strictly negative real part.
This yields the delay-independent asymptotic stability of \eqref{equa}.\\

Conversely, assume that \eqref{delt2} is not fulfilled.
Then,
\[
\Delta_{i'} - \sum_{i\in\cI\setminus\{ i'\}} \Delta_i
\]
admits a real positive eigenvalue.
In these conditions, system \eqref{equa} is not delay-independently stable, see \cite{HER84,HAL85}.
\end{proof}

\section{Examples}
\label{sec:examples}

\begin{figure}[t]
     \centerline{
       \includegraphics[width=4cm]{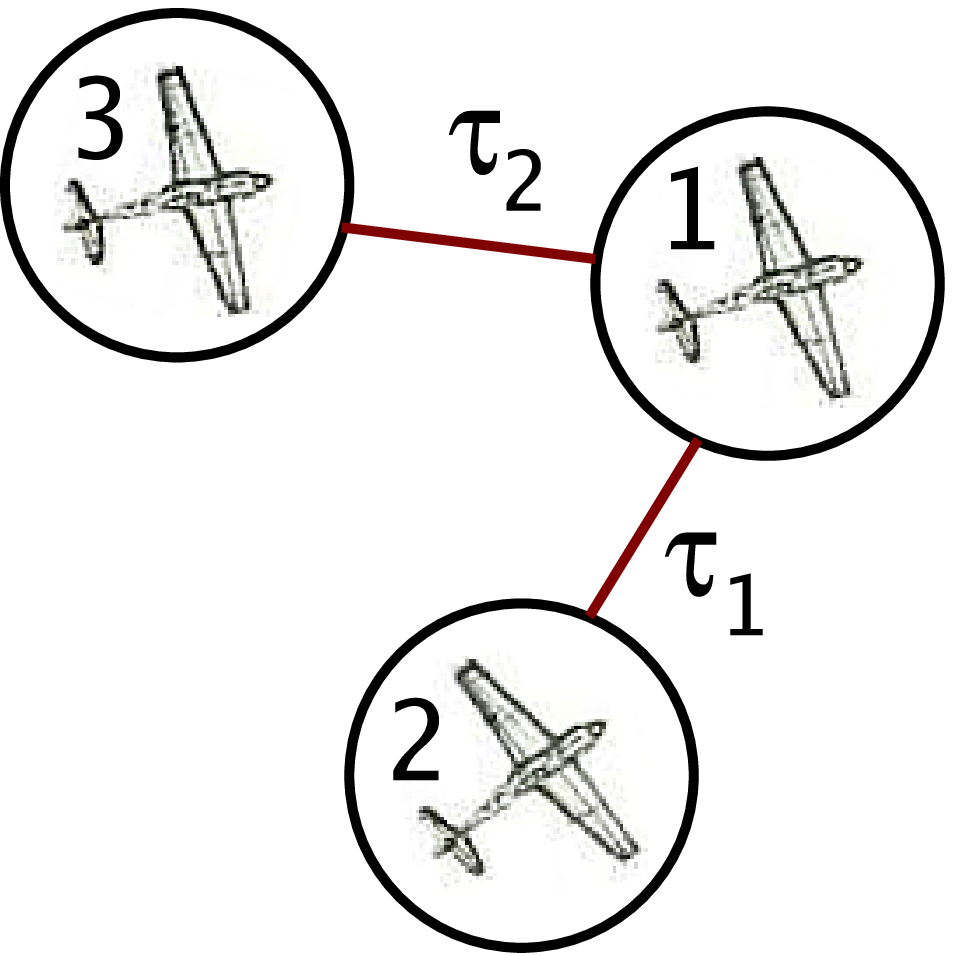}
}

     \caption{\protect{\small The multi-agent system, with the
         communication delays, used in Section~\ref{sec:examples}.}}
     \label{fig:graph}
  \end{figure}

We stress once more that the results in Sections \ref{se2} and
\ref{se3}
characterize robustness of average consensus, i.e. average consensus
{\em for any value of the delays less or equal to $\btau$}.
In order to illustrate this concept, we consider the network of three agents whose
communication graph $G$ is represented in Figure~\ref{fig:graph}. 

We assume that $v(x,t)\in\Rset^2$, that the weights $\omega(x,y)=1\Leftrightarrow x\sim y$
are used and that the delays
$\tau_i>0$, $i=1,2$ are constant in time. Moreover, the agents evolve
according to the PdE \eqref{eq:f_dynamics} starting from the initial
conditions
\begin{equation*}
  \tilde v(1,t)=[2,~2]', \quad
 \tilde v(2,t)=[2,~-2]', \quad
 \tilde v(3,t)=[1,~3]'
\end{equation*}
 where $t\in [-\max\{\tau_1,\tau_2\},0]$.
The average velocity at time $t=0$ is $\bv=[\frac{5}{3},~1]'$.

From Theorem \ref{th:eiglapl}, the eigenvalues of $\lapl$ on $H^1$
are the non null eigenvalues of $\MM(\lapl)$ (up to their multiplicity).
In our case, one gets $\|\lapl\|=3$, and the bound \eqref{cons} 
is equal to $\pi/6\simeq 0.524$.

\begin{figure}[t]
     \centerline{
     \begin{tabular}[t]{cc}
       \subfigure[$\tau_1=\tau_2=0.51$]{
       \includegraphics[width=7cm]{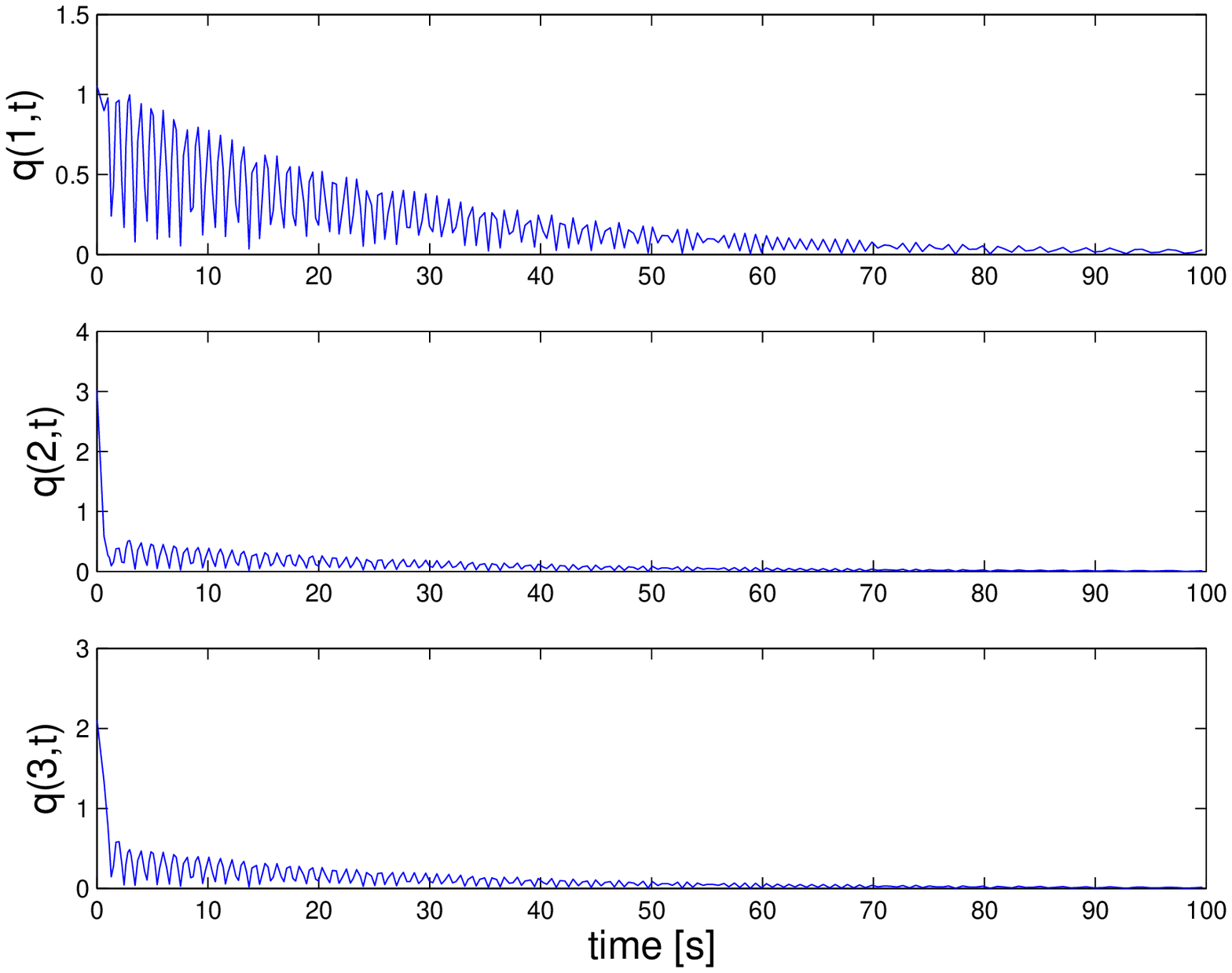}\label{fig:err-a}}&
       \subfigure[$\tau_1=\tau_2=0.53$]{
       \includegraphics[width=7cm]{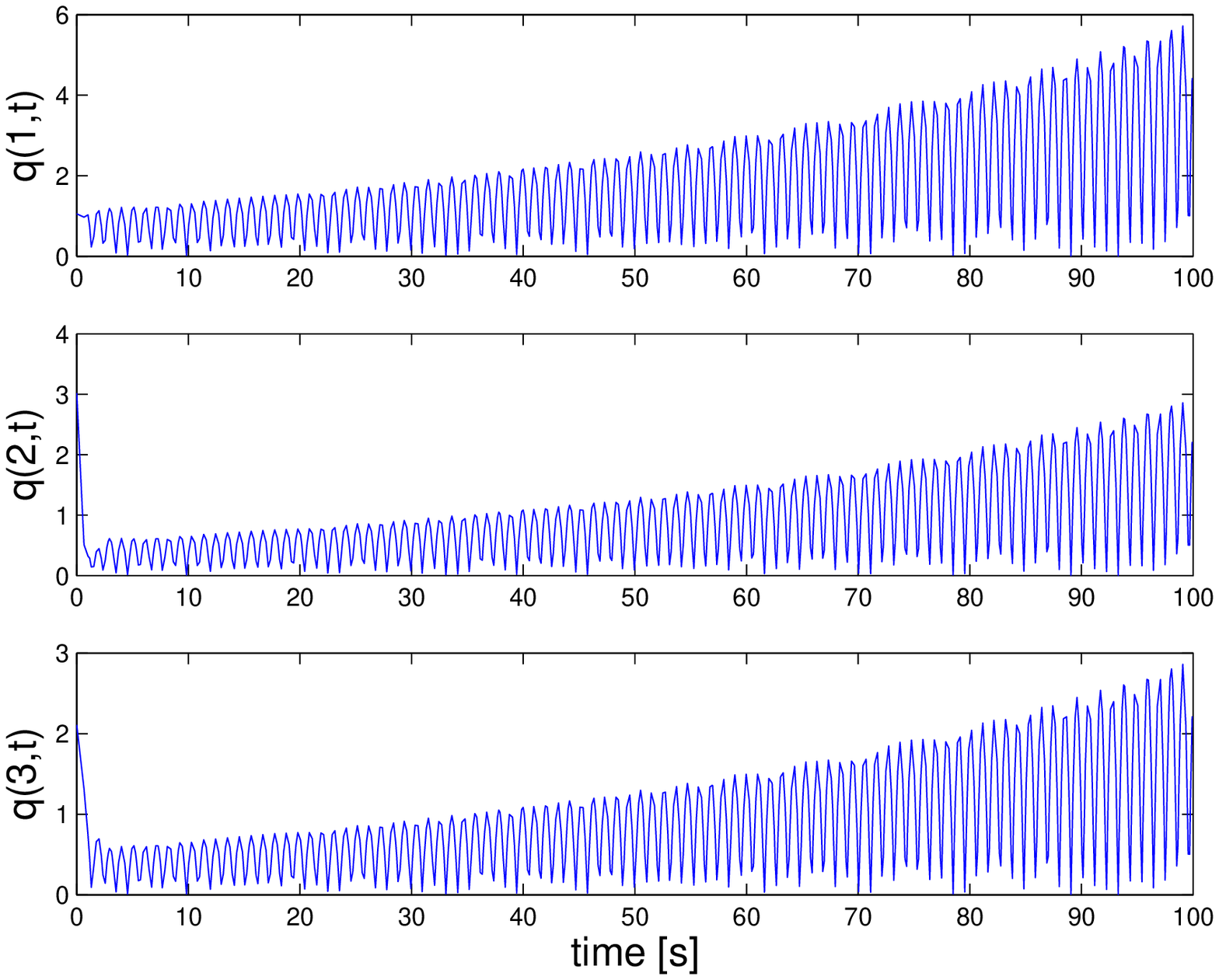}\label{fig:err-b}}
     \end{tabular}
     }
     \caption{\protect{\small Time evolution of $q(x,t)=\|v(x,t)-\bv\| $ for
         the multi-agent system in Figure~\ref{fig:graph}.}}
     \label{fig:err}
  \end{figure}

\begin{figure}[t]
     \centerline{
       \includegraphics[width=7cm]{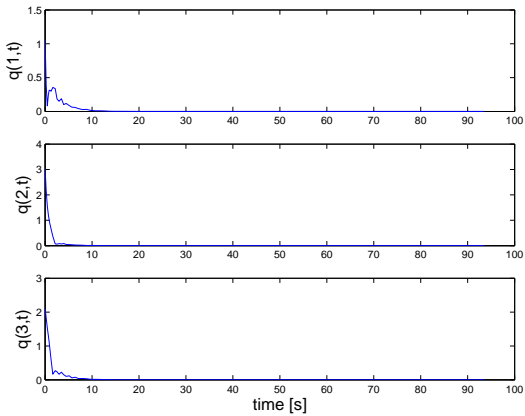}
}

     \caption{\protect{\small Time evolution of $q(x,t)=\|v(x,t)-\bv\|  $ for
         the multi-agent system in Figure~\ref{fig:graph}, with $\tau_1=0.1$ and $\tau_2=0.7$.}}
     \label{fig:err-c}
  \end{figure}

In the first experiment, we choose the delays $\tau_1=\tau_2=0.51$ that are
slightly below $\btau$. Then, Theorem~\ref{th2} guarantees average consensus
and such a result can be verified from Figure~\ref{fig:err-a}, where
the evolution of $\|v(x,t)-\bv\| $, $x\in\{1,2,3\}$ is represented.
In the second experiment, we use $\tau_1=\tau_2=0.53$, so having
$\tau_1=\tau_2>\btau$. As predicted by Theorem~\ref{th2}, the dynamics of
$v_1$ becomes unstable and average consensus cannot be achieved. This can be
clearly seen in Figure~\ref{fig:err-b}.
Finally, we choose $\tau_1=0.1$ and $\tau_2=0.7$. In this case,
$\tau_2$ violates the bound of Theorem~\ref{th2}. However, $\tau_1<\btau$
and Theorem~\ref{th2} cannot be used for checking the average consensus property. In this
case,
the achievement of average consensus can be verified by simulation, as
shown in Figure~\ref{fig:err-c}. 


\subsection{Two extremal cases}
\label{se5}

We now consider the case of complete and loop-shaped graphs with uniform weights.
In Sections \ref{subsec:complete} and
\ref{subsec:loop} it is shown that
the bounds of Table~\ref{ta1} can be computed in closed-form as
a function of $N$. 


\subsubsection{Complete graph}
\label{subsec:complete}

We consider a complete graph $G$ characterized by
$\EE=\NN\times\NN$ and $\omega(x,y)=\delta>0$.
We assume that a single delay $\tau_{x,y}$ is associated to each edge. Since delays are
symmetric, then $\DD=\{\tau_{x,y},(x,y)\in\II\}$ where
  $\II=\{(x,y)\in G\times G: x<y\}$ and the cardinality of $r$ of $\II$, defined in \eqref{eq:dset}, is 
$\frac{N(N-1)}{2}$. Note that the only difference
  with respect to \eqref{eq:dset} is that delays are parametrized by 
two indexes. One also has that
  $G_{x,y}=(\NN,T^{-1}(\tau_{x,y})=\{(x,y),(y,x)\}$ and, by direct
  computation, 
\begin{equation}
\label{eq:lapl_xy}
(\lapl_{x,y}f)(\xi)=
\begin{cases}
0 &\mbox{ if } \xi\not \in G_{x,y} \\
\delta(f(y)-f(x)) &\mbox{ if } \xi=x \\
\delta(f(x)-f(y)) &\mbox{ if } \xi=y \\
\end{cases}
 \end{equation} 
The matrix representation of $\lapl_{x,y}$ is given by 
\begin{equation}
\MM(\lapl_{x,y})=-\delta e_{x,y}e_{x,y}^T \otimes I_d \label{eq:mlapl_xy}
\end{equation}
where $e_{x,y}$ is defined
as the vector in $\Rset^N$ with zero entries, except the $x$-th one, equal to $1$,
and the $y$-th one, equal to $-1$.
Note that, from \eqref{eq:lapl-split}, we have 
\begin{equation}
  \label{eq:lapl-split2}
  \lapl=\sum_{x<y}\lapl_{x,y}
\end{equation}
For computing the bound in \eqref{eq:bound-vary-nonuniform}, one needs
to evaluate $\|\lapl_{x',y'}\lapl_{x,y}\|$ for $x<y$ and $x'<y'$. 
Note that
$(\lapl_{x',y'}\lapl_{x,y}f)(\xi)$ is possibly non null only if  the indexes
$(x',y')$ and $(x,y)$ have both elements in common
(i.e. $(x',y')=(x,y)$) 
or if they have just one element in
common. The latter occours if and only if
{\small
\begin{equation}
x<y, x'<y' \mbox{ and }
\left(
 x=x', y\neq y'  \mbox { or }  x\neq x', y=y'  \mbox { or } 
x=y', y\neq x'   \mbox { or }  x\neq y', y=x' 
\right) \label{eq:linked}
\end{equation}
}
The next lemma provides closed-form expression for the quantities
appearing in Table~\ref{ta1}.
\begin{lemma}
\label{le2}
For a complete graph $G$ with uniform  weights $\omega(x,y)=\delta>0$
it holds
\begin{subequations}
\begin{gather}
\label{forma}
\left\|\lapl\right\| =
\left\|\lapl^{-1}
\right\|^{-1} = N\delta, \\
\label{formb}
\left\|
e_{x,y}e_{x,y}^Te_{x',y'}e_{x',y'}^T
\right\| =
\begin{cases}
4 & \text{ if } (x',y')=(x,y)\\
1 & \text{ if } (x',y') \mbox{ and } (x,y) \mbox{ verify }\eqref{eq:linked}, \\
0 & \text{ if } x\neq x' \mbox{ and } y\neq y',
\end{cases}\\
\label{formc}
\sum_{x<y, x'<y'} \left\|\lapl_{x',y'}\lapl_{x,y}
\right\| = \delta^2 N^2(N-1) \ .
\end{gather}
\end{subequations}
\end{lemma}

\begin{proof}
In \cite[pag. 269]{B98} it is shown that the non-zero eigenvalues of
$\MM(\lapl)$ are all equal to $-N$, when
$\delta=1$. In view of Theorem~\ref{th:eiglapl}, these are the
eigenvalues of $\lapl$ on $H^1$ and \eqref{forma}
follows.
A simple proof of the same fact is given next. From the graph completeness,
one has $\lapl f(x)=-\delta N
f(x) +\delta \int_G f(y)~dy$. Then, $\lapl f=\lambda f$, $f\in H^1$
results in
\[
-\delta N
f(x) +\delta \int_G f(y)~dy=\lambda f(x)
\]
and since the average of $f$ is zero,
\[
-\delta N
f(x) =\lambda f(x)
\]
so proving that the eigenvalues of $\lapl$ on $H^1$ are all equal to
$-N\delta$.


Formula \eqref{formb} is obtained directly.
When $(x',y')=(x,y)$  we have that
$(e_{x,y}e_{x,y}^T)^2$ has all zero entries except the
elements $(x,x)$, $(y,y)$ (that are equal to $2$) and $(x,y)$,
$(y,x)$ (that are equal to $-2$). The minimal eigenvalue of $(e_{x,y}e_{x,y}^T)^2$
is always $-4$.  
When $(x',y')$ and $(x,y)$ verify \eqref{eq:linked}, one can show that
up to a permutation of the agent indexes, which does not affect the results, it holds

\[
\|e_{x,y}e_{x,y}^Te_{x',y'}e_{x',y'}^T\|=\left\|
\begin{pmatrix}
1 & -1 & 0 \\ -1 & 1 & 0\\ 0 & 0 & 0
\end{pmatrix}
\begin{pmatrix}
0 & 0 & 0\\ 0 & 1 & -1 \\ 0 & -1 & 1
\end{pmatrix}
\right\|
=
\left\|
\begin{pmatrix}
0 & 1 & -1 \\ 0 & -1 & 1\\0 & 0 & 0
\end{pmatrix}
\right\|
= 1\ .
\]

Finally, from \eqref{formb}, one gets
\begin{equation}
\label{eq:sums}
\sum_{x<y, x'<y'} \left\|
e_{x,y}e_{x,y}^Te_{x',y'}e_{x',y'}^T
\right\|
= 4\sum_{x=x'<y=y'}1+1\sum_{\scriptsize \begin{array}{c} (x,y),(x',y') \\
\mbox{\scriptsize{verifying\eqref{eq:linked}} }\end{array}}1
\end{equation}
The first sum in the previous formula is equal to $\frac{N(N-1)}{2}$.
As for the second sum, for a fixed pair $(x,y)$, $1\leq x<y\leq N$  one can show that
the number of pairs $(x',y')$ verifying $1\leq x'<y'\leq N$, $x<y$, $x'<y'$ and
\begin{eqnarray}
x=x,\ y\neq y'&\mbox{  is  } &N-x-1 \label{eq:sl1} \\
x\neq x,\ y= y'&\mbox{  is  } &y-2 \label{eq:sl2} \\
x=y',\ y\neq x'&\mbox{  is  } &x-1 \label{eq:sl3} \\  
x\neq y',\ y= x'&\mbox{  is  } &N-y \label{eq:sl4} 
\end{eqnarray}
Then,
\begin{equation}
  \label{eq:count}
  \sum_{\scriptsize \begin{array}{c} (x,y),(x',y') \\
\mbox{\scriptsize{verifying\eqref{eq:linked}} }\end{array}}1=\frac{N(N-1)}{2}(N-x-1+y-2+x-1+N-y)=N(N-1)(N-2)
\end{equation}
and one finally gets
\begin{equation}
 \label{eq:sums2}
\sum_{x<y, x'<y'} \left\|
e_{x,y}e_{x,y}^Te_{x',y'}e_{x',y'}^T
\right\|=\frac{N(N-1)}{2}+N(N-1)(N-2)=N^2(N-1)
\end{equation}
From \eqref{eq:mlapl_xy}, one has that 
$\MM(\lapl_{x,y}\lapl_{x',y'})=\delta^2e_{x,y}e_{x,y}^Te_{x',y'}e_{x',y'}^T$.
and formula
\eqref{formc} follows.

\end{proof}

\subsubsection{Loop-shaped graph}
\label{subsec:loop}

We now analyze the case of a loop-shaped graph with uniform weights 
where each agent
exchanges information only with two other ones. More precisely, the set of
edges is given by $\EE=\{(1,2),(2,3),\ldots,(N-1,N),(N,1)\}$ and 
$\omega(x,y)=\delta>0$ if $(x,y)\in\EE$.
As in Section~\ref{subsec:complete}, we associate a different delay to
each edge. Then, $\DD=\{\tau_{x,y},(x,y)\in \II\}$, where
$\II=\EE$. The operators $\lapl_{x,y}$ are defined as in
\eqref{eq:lapl_xy} and their matrix representation is given in 
\eqref{eq:mlapl_xy}.
The matrix representation of $\lapl$ is given by
\begin{equation}
\MM(\lapl)=  \delta \left(
-2I_N+P_N+P_N^T
\right)\otimes I_d \label{eq:mlapl2}
\end{equation}.  
where $P_N\in\Rset^{N\times N}$ is the permutation matrix 
\begin{equation}
  \label{eq:perm}
\small
  P_N=\begin{pmatrix}
0 & 0 & \cdots & 0 & 1 \\
1 & 0 & \cdots & 0 & 0 \\
0 & 1 & \cdots & 0 & 0 \\
\vdots & \vdots  & \ddots & \vdots & \vdots   \\
0 & 0 & \cdots & 1 & 0
\end{pmatrix}
\end{equation}

The decomposition \eqref{eq:lapl-split} becomes
\[
\lapl = - \delta
\sum_{(x,y)\in\cI} e_{x,y}e_{x,y}^T \otimes I_d\ ,
\]
where the definition of $e_{x,y}$, has been given in Section \ref{subsec:complete}.

The counterpart of Lemma \ref{le2} is stated next.

\begin{lemma}
\label{le3}
For a loop-shaped  graph $G$ with uniform  weights $\omega(x,y)=\delta>0$
it holds
\begin{subequations}
\begin{gather}
\label{formd}
\left\|
\lapl
\right\| = 4\sin^2 \left(\left\lfloor
\frac{N}{2}
\right\rfloor
\frac{\pi}{N}
\right), \quad
\left\|\left(
\lapl
\right)^{-1} \right\|^{-1} = 4\delta \sin^2\left(\frac{\pi}{N}\right),\\
\label{frome}
\sum_{(x,y),(x',y')\in \II} \left\|\lapl_{x',y'}\lapl_{x,y}
\right\| = 6\delta^2N.
\end{gather}
\end{subequations}
\end{lemma}

\begin{proof}
As is well-known, $P_N^T=P_N^{-1}$, and the eigenvalues $\lambda_{P_N}$ of
$P_N$ are all the $N$-th roots of the unit, i.e. $\lambda_{P_N}=\exp(\frac{2\pi
  k}{N})$, $k=1,\ldots,N$.
An eigenvalue $\lambda$ of $\MM(\lapl)$ fulfills the equation 
\[
(P_N+P_N^{-1})\psi=\tilde \lambda \psi,\quad \tilde \lambda=\lambda+2
\]
from which $\tilde \lambda=\lambda_{P_N}+\lambda_{P_N}^{-1}=2\cos\left(\frac{2\pi
    k}{N}\right)$. Then,
$\lambda=2\cos\left(\frac{2\pi k}{N}\right)-2=-4\sin^2 \frac{k\pi}{N}$. From
Theorem \ref{th:eiglapl}, a nonzero $\lambda$ is an eigenvalue of
$\lapl$ on $H^1$.
This allows to deduce the identities in \eqref{formd}.

From identity \eqref{formb}, it turns out that
\begin{equation}
\label{eq:sums3}
\sum_{(x,y)\in\II,(x',y')\in\II} \left\|
e_{x,y}e_{x,y}^Te_{x',y'}e_{x',y'}^T
\right\|
= 4\sum_{(x,y)=(x'y')\in\II}1+1\sum_{\scriptsize \begin{array}{c} (x,y),(x',y') \\
\mbox{\scriptsize{verifying \eqref{eq:linked}} }\end{array}}1
\end{equation}
The first sum in \eqref{eq:sums3} is equal to $N$.
By direct calculation, the number of  pairs $(x,y)\in\II$ that verify \eqref{eq:linked} for a given pair $(x',y')\in\II$ is $2$. Therefore,
since the number of possible pairs $(x',y')\in \II$ is
$N$, one has that the total number of pairs $(x',y')$
and $(x,y)$ that verify \eqref{eq:linked} is $2N$.
Then, 
one obtains 
\[
\sum_{(x,y)\in\II,(x',y')\in\II} \left\|
e_{x,y}e_{x,y}^Te_{x',y'}e_{x',y'}^T
\right\|
=4N+2N=6N
\]
From \eqref{eq:mlapl_xy}, one has that 
$\MM(\lapl_{x,y}\lapl_{x',y'})=\delta^2e_{x,y}e_{x,y}^Te_{x',y'}e_{x',y'}^T$.
and formula \eqref{frome} follows.
\end{proof}



The delay margins obtained by applying the results summarized in Table \ref{ta1}
are given now in Tables \ref{ta2} and \ref{ta3}, for each of the two networks.

\begin{table}[h]
\begin{center}
\begin{tabular}{c||c|c|}
$\btau$
& Uniform delays & Non-uniform delays \\
\hline
& & \\
Time-invariant delays & $\frac{\pi}{2\delta N}$
& $\frac{\pi}{2\delta N}$\\
& & \\
\hline
& & \\
Time-varying delays & $\frac{3}{2\delta N}$
& $\frac{1}{\delta N (N-1)}$\\
& &\\
\hline
\end{tabular}
\caption{ Complete graph: Bounds on the worst-case stabilizing delays}
\label{ta2}
\end{center}
\end{table}

\begin{table}[h]
\begin{center}
\begin{tabular}{c||c|c|}
$\btau$
& Uniform delays & Non-uniform delays \\
\hline
& & \\
Time-invariant delays & $\frac{\pi }{8\delta \sin^2 \left(\left\lfloor
\frac{N}{2}
\right\rfloor
\frac{\pi}{N}
\right)}$
& $\frac{\pi }{8\delta \sin^2 \left(\left\lfloor
\frac{N}{2}
\right\rfloor
\frac{\pi}{N}
\right)}$\\
& & \\
\hline
& & \\
Time-varying delays &  $\frac{3 }{8\delta \sin^2 \left(\left\lfloor
\frac{N}{2}
\right\rfloor
\frac{\pi}{N}
\right)}$
& $\frac{2\sin^2 \frac{\pi}{N}}{3\delta N}$\\
& &\\
\hline
\end{tabular}
\caption{Loop-shaped graph: Bounds on the worst-case stabilizing
  delays.}
\label{ta3}
\end{center}
\end{table}

\section{Conclusions}

We provided convergence analysis of an average consensus protocol
for undirected networks of dynamic agents having communication delays.
We considered constant or time-varying delays, uniformly or non uniformly distributed in the network.
Sufficient conditions (also necessary in most cases) for existence of average consensus
under bounded, but otherwise unknown, communication delays, have been given. 
Simulations have been provided that demonstrate adequation with the bounds computed analytically.

\bibliographystyle{plain}
\bibliography{flock}








\end{document}